\documentclass[12pt,twoside,leqno]{article}
\usepackage{latexsym}
\def\implies{\; \Longrightarrow \;}

\begin{document}
\begin{center}
{\bf ON THE SEMIGROUP OF ORDER-DECREASING PARTIAL ISOMETRIES OF A FINITE CHAIN}\\[4mm]
\textbf{R. Kehinde, S. O. Makanjuola and A. Umar}\\

\end{center}

\newtheorem{theorem}{{\bf Theorem}}[section]
\newtheorem{prop}[theorem]{{\bf Proposition}}
\newtheorem{lemma}[theorem]{{\bf Lemma}}
\newtheorem{corollary}[theorem]{{\bf Corollary}}
\newtheorem{remark}[theorem]{{\bf Remark}}
\newtheorem{conj}[theorem]{{\bf Conjecture}}
\newcommand{\pf}{\smallskip\noindent {\em Proof.}\ \  }
\newcommand{\qed}{\hfill $\Box$\medskip}

\newcommand{\inv}{^{-1}}

\begin{abstract}
Let ${\cal I}_n$ be the symmetric inverse semigroup on $X_n = \{1,
2, \cdots , n\}$ and let ${\cal DDP}_n$ and ${\cal ODDP}_n$ be its
subsemigroups of order-decreasing partial isometries and of
order-preserving order-decreasing partial isometries of $X_n$,
respectively. In this paper we investigate the cycle structure of
order-decreasing partial isometry and characterize the Green's
relations on ${\cal DDP}_n$ and ${\cal ODDP}_n$. We show that ${\cal
ODDP}_n$ is a $0-E-unitary$ ample semigroup. We also investigate the
cardinalities of some equivalences on ${\cal DDP}_n$ and ${\cal
ODDP}_n$ which lead naturally to obtaining the order of the
semigroups.\footnote{\textit{Key Words}: partial one-one
transformation, partial isometries, height, right (left) waist,
right (left) shoulder and fix of a transformation, idempotents and
nilpotents.} \footnote{This work was carried out when the first
named author was visiting Sultan Qaboos University for a 3-month
research visit in Fall 2010.}\end{abstract}

\textit{MSC2010}: 20M18, 20M20, 05A10, 05A15.

\section{Introduction and Preliminaries}

Let $X_n=\{1,2, \cdots, n\}$ and ${\cal I}_n$ be the partial
one-to-one transformation semigroup on $X_n$ under composition of
mappings. Then ${\cal I}_n$ is an {\em inverse} semigroup (that is,
for all $\alpha \in {\cal I}_n$ there exists a unique $\alpha' \in
{\cal I}_n$ such that $\alpha = \alpha\alpha'\alpha$ and $\alpha' =
\alpha'\alpha\alpha'$). The importance of ${\cal I}_n$ (more
commonly known as the symmetric inverse semigroup or monoid) to
inverse semigroup theory may be likened to that of the symmetric
group ${\cal S}_n$ to group theory. Every finite inverse semigroup
$S$ is embeddable in ${\cal I}_n$, the analogue of Cayley's theorem
for finite groups, and to the regular representation of finite
semigroups. Thus, just as the study of symmetric, alternating and
dihedral groups has made a significant contribution to group theory,
so has the study of various subsemigroups of ${\cal I}_n$, see for
example \cite{Bor, Fer1, Fer2, Gan, Gar, Uma1}.

\noindent A transformation $\alpha \in {\cal I}_n$ is said to be a
{\em partial isometry} if (for all $x,y\in {Dom\,\alpha}$) $\mid x-y
\mid = \mid x\alpha -y\alpha\mid$; {\em order-preserving
(order-reversing)} if (for all $x,y \in \ {Dom\, \alpha})\ x \leq y
\implies x\alpha \leq y\alpha \ (x\alpha \geq y\alpha)$; and, is
said to be {\em order-decreasing} if (for all $x\in {Dom\,
\alpha})\, x\alpha \leq x$. Semigroups of partial isometries on more
restrictive but richer mathematical structures have been studied
\cite{Bra, Wal}. Recently, the authors in \cite{Keh} studied the
semigroup of partial isometries of a finite chain, ${\cal DP}_n$ and
its subsemigroup of order-preserving partial isometries ${\cal
ODP}_n$. Ealier, one of the authors studied the semigroup of partial
one-to-one order-decreasing(order-increasing) transformations of a
finite chain, ${\cal I}_n^-$ \cite{Uma1}. This paper investigates
the algebraic and combinatorial properties of ${\cal DDP}_n$ and
${\cal ODDP}_n$, the semigroups of order-decreasing partial
isometries and of order-preserving order-decreasing partial
isometries of an $n-$chain, respectively.

In this section we introduce basic terminologies and some
preliminary results concerning the cycle structure of a partial
order-decreasing isometry of $X_n$. In the next section, (Section 2)
we characterize the classical Green's relations and their starred
analogues, where we show that ${\cal ODDP}_n$ is a (nonregular)
0-E-unitary ample semigroup. We also show that certain Rees factor
semigroups of ${\cal ODDP}_n$ are 0-E-unitary and categorical ample
semigroups. In Section 3 we obtain the cardinalities of two
equivalences defined on ${\cal DDP}_n$ and ${\cal ODDP}_n$. These
equivalences lead to formulae for the order of ${\cal DDP}_n$ and
${\cal ODDP}_n$ as well as new triangles of numbers not yet recorded
in \cite{Slo}.

For standard concepts in semigroup and symmetric inverse semigroup
theory, see for example \cite{How, Lim, Law1}. In particular E(S)
denotes the set of idempotents of S. Let
\begin{eqnarray} \label{eqn1.1} {\cal DDP}_n= \{\alpha \in {\cal DP}_n:
(\forall \ x\in{Dom \ \alpha}) \ x\alpha \leq x\}.
\end{eqnarray}
\noindent be the subsemigroup of ${\cal I}_n$ consisting of all
order-decreasing partial isometries of $X_n$. Also let
\begin{eqnarray} \label{eqn1.2} {\cal ODDP}_n= \{\alpha \in {\cal DDP}_n:
(\forall \ x,y\in {Dom \ \alpha})\ x\leq y \Longrightarrow
x\alpha\leq y\alpha\}
\end{eqnarray}
\noindent be the subsemigroup of ${\cal DDP}_n$ consisting of all
order-preserving order-decreasing partial isometries of $X_n$. Then
we have the following result.

\begin{lemma} \label{lem1.1} ${\cal DDP}_n$ and ${\cal ODDP}_n$ are
subsemigroups of ${\cal I}_n$.
\end{lemma}

\begin{remark}\label{rem1} ${\cal DDP}_n={\cal DP}_n \cap{\cal I}_n^-$ and
${\cal ODDP}_n={\cal ODP}_n \cap{\cal I}_n^-$, where ${\cal I}_n^-$
is a semigroup of partial one-to-one order-decreasing
transformations of $X_n$.
\end{remark}

As in \cite{Keh}, we prove a sequence of lemmas that help us
understand the cycle structure of order-decreasing partial
isometries. These lemmas also seem to be useful in investigating the
combinatorial questions in Section 3. First, let $\alpha$ be in
${\cal I}_n$. Then the {\em height} of $\alpha$ is $h(\alpha)=
\mid{Im\,\alpha}\mid$, the {\em right [left] waist} of $\alpha$ is
$w^+(\alpha) = max({Im \, \alpha})\, [w^-(\alpha) = min({Im
\,\alpha})]$, the {\em right [left] shoulder} of $\alpha$ is
$\varpi^+(\alpha) = max({Dom \, \alpha})$\, [$\varpi^-(\alpha) =
min({Dom \, \alpha})]$, and {\em fix} of $\alpha$ is denoted by
$f(\alpha)$, and defined by $f(\alpha)=|F(\alpha)|$, where
$$F(\alpha) = \{x \in X_n: x\alpha = x\}.$$

\begin{lemma} \cite[Lemma 1.2]{Keh} \label{lem1.2} Let $\alpha\in {\cal DP}_n$
be such that $h(\alpha)=p$. Then $f(\alpha)=0 \,or\, 1\, or\, p$.
\end{lemma}

\begin{corollary} \cite[Corollary 1.3]{Keh} Let $\alpha\in {\cal DP}_n$.
If $f(\alpha)=p>1$ then $f(\alpha)=h(\alpha)$. Equivalently, if
$f(\alpha)>1$ then $\alpha$ is an idempotent.
\end{corollary}

\begin{lemma} \label{lem1.3} Let $\alpha\in {\cal DDP}_n$.
If $i\in F(\alpha)$ ($1\leq i\leq n)$ then for all $x\in {Dom \
\alpha}$, such that $x \ < i$ we have $x\alpha =x$.
\end{lemma}

\pf Note that for all $x\in {Dom\,\alpha}$ we have $x\alpha\leq x <
i$ and so $i-x=\\\mid i\alpha - x\alpha\mid=\mid i -
x\alpha\mid=i-x\alpha \Longrightarrow x=x\alpha.$\qed

\begin{corollary} \label{cor1.4} Let $\alpha\in {\cal DDP}_n$.
If $F(\alpha)=\{i\}$ then ${Dom\,\alpha}\subseteq \\ \{i, i+1,
\cdots ,n\}$.
\end{corollary}

\begin{lemma} \cite[Lemma 1.4]{Keh}\label{lem1.5} Let $\alpha\in {\cal DP}_n$. If
$1\in F(\alpha)$ or $n\in F(\alpha)$ then for all $x\in {Dom
\alpha}$, we have $x\alpha =x$. Equivalently, if $1\in F(\alpha)$ or
$n\in F(\alpha)$ then $\alpha$ is a partial identity.
\end{lemma}

\begin{lemma} \cite[Lemma 1.5]{Keh} \label{lem1.6} Let $\alpha\in {\cal ODP}_n$ and
 $n\in {Dom\, \alpha}\cap{Im\, \alpha}$. Then $n\alpha =n$.
\end{lemma}

\begin{lemma} \cite[Lemma 1.6]{Keh} \label{lem1.7} Let $\alpha\in {\cal ODP}_n$
and $f(\alpha)\geq 1$. Then $\alpha$ is an idempotent.
\end{lemma}

\begin{lemma} \label{lem1.8} Let $\alpha\in {\cal ODDP}_n$. Then
$x-x\alpha=y-y\alpha$ for all $x,y\in {Dom\,\alpha}$.
\end{lemma}

\pf let $x,y\in {Dom\,\alpha}$ be such that $x>y$. Then by the
order-preserving and isometry properties we see that $\mid
x-y\mid=\mid x\alpha- y\alpha\mid \Longrightarrow
x-y=x\alpha-y\alpha \Longrightarrow x-x\alpha=y-y\alpha.$ The case
$x<y$ is similar.\qed

\section{Green's relations and their starred analogues}

For the definitions of Green's relations we refer the reader to
Howie \cite[Chapter 2]{How2}. First we have

\begin{theorem} \label{thrm2.1} Let ${\cal DDP}_n$ and ${\cal ODDP}_n$ be as defined in
(\ref{eqn1.1}) and (\ref{eqn1.2}) respectively. Then ${\cal DDP}_n$
and ${\cal ODDP}_n$ are\ ${\cal J}$-trivial.
 \end{theorem}

\pf It follows from \cite[Lemma 2.2]{Uma1} and Remark
\ref{rem1}.\qed

Now since ${\cal ODDP}_n$ contains some nonidempotent elements:
$$\pmatrix{x\cr y}(x > y)\,$$ it follows immediately that

\begin{corollary} \label{cor1} For $n>1$, ${\cal DDP}_n$ and ${\cal ODDP}_n$
are non-regular semigroups.
\end{corollary}

On the semigroup $S$ the relation ${\cal L^*}({\cal R^*})$ is
defined by the rule that $(a,b)\in {\cal L^*}({\cal R^*})$ if and
only if the elements $a,b$ are related by the Green's relation
${\cal L}({\cal R})$ in some oversemigroup of $S$. The join of the
equivalences ${\cal L^*}$ and ${\cal R^*}$ is denoted by ${\cal
D^*}$ and their intersection by ${\cal H^*}$. For the definition of
the starred analogue of the Green's relation ${\cal J}$, see
\cite{Fou2} or \cite{Uma1}.

A semigroup $S$ in which each ${\cal L^*}$-class and each ${\cal
R^*}$-class contains an idempotent is called $abundant$ \cite{Fou2}.

By \cite[Lemma1.6]{Qal} and \cite[Proposition 2.4.2 \& Ex.
5.11.2]{How2} we deduce the following lemma.

\begin{lemma} Let $\alpha ,\beta \in {\cal DDP}_n.$ \ Then
\begin{itemize}
\item[(1)] $\alpha \leq _{{\cal R^*}}\beta $ \,  if and only if\,
${Dom\,\alpha}\subseteq {Dom\,\beta}$;
\item[(2)] $\alpha \leq _{{\cal L^*}}\beta $\,   if and only if  $Im\,\alpha
\subseteq Im\,\beta $;
\item[(3)] $\alpha \leq _{{\cal H^*}}\beta $ \,  if and only if\,
${Dom\,\alpha}\subseteq {Dom\,\beta}$ and $Im\,\alpha \subseteq
Im\,\beta $.
\end{itemize}
\end{lemma}

\pf It is enough to observe that ${\cal ODDP}_n$ and ${\cal DDP}_n$
are full subsemigroups of ${\cal I}_n$ in the sense that $E({\cal
ODDP}_n)$=$E({\cal DDP}_n)$=$E({\cal I}_n)$.\qed

An abundant semigroup $S$ in which $E(S)$ is a semilattice is called
$adequate$ \cite{Fou1}. Of course inverse semigroups are adequate
since in this case ${\cal L^*}={\cal L}$ and ${\cal R^*}={\cal R}$.

As in \cite{Fou1}, for an element $a$ of an adequate semigroup $S$,
the (unique) idempotent in the ${\cal L^*}$-class(${\cal
R^*}$-class) containing $a$ will be denoted by $a^*(a^+)$. An
adequate semigroup $S$ is said to be {\em ample} if $ea$=$a(ea)^*$
and $ae$=$(ae)^+a$ for all elements $a$ in $S$ and all idempotents
$e$ in $S$. Ample semigroups were known as $type A$ semigroups.

\begin{theorem} Let ${\cal DDP}_n$ and ${\cal ODDP}_n$ be as defined in (\ref{eqn1.1})
and (\ref{eqn1.2}) respectively. Then ${\cal DDP}_n$ and ${\cal
ODDP}_n$ are non-regular ample semigroups.
\end{theorem}

\pf The proofs are similar to that of \cite[theorem 2.6]{Uma1}.\qed

\begin{theorem} Let $S = {\cal ODDP}_n$ be as defined in (\ref{eqn1.2}).
Then $\alpha \leq _{{\cal D^*}}\beta $ if and only if there exists
an order-preserving isometry $\theta : {Dom\,\alpha}\rightarrow
{Im\,\beta}.$
\end{theorem}

\noindent Let $E'=E \setminus {0}$. A semigroup S is said to be
$0-E-unitary$ if $(\forall e\in E')(\forall s\in S)\,\,\, es\in E'
\Longrightarrow s\in E'$. The structure theorem for 0-E-unitary
inverse semigroup was given by Lawson \cite{Law2}, see also Szendrei
\cite{Sze} and Gomes and Howie \cite{Gom}.

\begin{theorem} \label{thrm3.3} ${\cal ODDP}_n$ is a $0-E-unitary$
ample subsemigroup of ${\cal I}_n$.
\end{theorem}

\pf It follows from \cite[Theorem 2.4]{Keh}.\qed

\begin{remark} Note that ${\cal DDP}_n$ is not 0-E-unitary:
$$\pmatrix{1&2\cr 1&2}\pmatrix{2&3\cr
2&1}=\pmatrix{2\cr 2}\in E({\cal
DDP}_n)\,\,\mbox{but}\,\,\pmatrix{2&3\cr 2&1}\notin E({\cal
DDP}_n).$$
\end{remark}

\noindent For natural numbers $n,p$ with $n\geq p\geq 0$, let
\begin{eqnarray}  L(n,p)=\{\alpha\in {\cal ODDP}_n:h(\alpha)\leq p\}
\end{eqnarray}
\noindent be a two-sided ideal of ${\cal ODDP}_n$, and for $p>0$,
let
\begin{eqnarray} \label{eqn2.1} Q(n,p)= L(n,p)/L(n,p-1)
\end{eqnarray}
\noindent be its Rees quotient semigroup. Then $Q(n,p)$ is a
0-E-unitary semigroup whose nonzero elements may be thought of as
the elements of ${\cal ODDP}_n$ of height $p$. The product of two
elements of $Q(n,p)$ is 0 whenever their product in ${\cal ODDP}_n$
is of height less than $p$.

\noindent A semigroup S is said to be $categorical$ \cite{Gom} if
$$(\forall a,b,c\in S),\,\, abc=0 \Longrightarrow ab=0\,\mbox{or}\,
bc=0$$.

\begin{theorem} \label{thrm3.4} Let $Q(n,p)$ be as defined in (\ref{eqn2.1}). Then
$Q(n,p)$ is a $0-E-unitary$ categorical semigroup.
\end{theorem}

\pf It follows from \cite[thrm2.6]{Keh}.\qed

\begin{remark} Note that ${\cal ODDP}_n$ is not categorical: $$\pmatrix{1&2\cr
1&2}\pmatrix{2&3\cr 2&3}\pmatrix{1&3\cr 1&3}= 0$$ \noindent but
$$\pmatrix{1&2\cr 1&2}\pmatrix{2&3\cr 2&3}= \pmatrix{2\cr 2}\neq 0
\,\,\,\mbox{and}\,\,\,\pmatrix{2&3\cr 2&3}\pmatrix{1&3\cr 1&3}=
\pmatrix{3\cr 3}\neq 0.$$
\end{remark}

\section{Combinatorial results}

For a nice survey article concerning combinatorial problems in the
symmetric inverse semigroup and some of its subsemigroups we refer
the reader to Umar \cite{Uma2}.

\noindent Now recall the definitions of {\em height} and {\em fix}
of $\alpha\in {\cal I}_n$ from the paragraph after Lemma
\ref{lem1.1}. As in Umar \cite{Uma2}, for natural numbers $n\geq
p\geq m\geq 0$ we define

\begin{eqnarray} \label{eqn3.1} F(n;p)= \mid\{\alpha \in S:
h(\alpha)=\mid {Im\,\alpha}\mid =p \}\mid,
\end{eqnarray}

\begin{eqnarray} \label{eqn3.2} F(n;m)= \mid\{\alpha \in S:
f(\alpha)=m \}\mid
\end{eqnarray}

\noindent where $S$ is any subsemigroup of ${\cal I}_n$. Also, let
$i=a_i=a$, for all $a \in \{p,m\}$, and $0\leq i\leq n$.

\begin{lemma}\label{lem3.1} Let $S={\cal ODDP}_n$. Then
$F(n;p_1)=F(n;1)=\pmatrix{n+1\cr 2}$ and $F(n;p_n)=F(n;n)=1$, for
all $n\geq 1$.
\end{lemma}

\pf Consider $\alpha$=$\pmatrix{x\cr x \alpha}$, where $x\geq x
\alpha$. If $x \alpha$=$i$ then $x \in \{i,i+1, \cdots, n\}$ and so
$x$ has $n-i+1$ degrees of freedom. Hence there are $\sum_{i=1}^{n}\
(n-i+1)=\frac{n(n+1)}{2}=\pmatrix{n+1\cr 2}$, order-decreasing
partial isometries of height 1. For the second statement, it is not
difficult to see that there is exactly one order-decreasing partial
isometry of height $n$: $\pmatrix{1&2&\dots&n\cr 1&2&\dots &n}$ (the
identity). \qed

\begin{lemma} \label{lem3.2} Let $S={\cal ODDP}_n$. Then
$F(n;p)=F(n-1;p-1)+F(n-1;p)$, for all $n\geq p\geq 2$.
\end{lemma}

\pf Let $\alpha \in {\cal ODDP}_n$ and $h(\alpha)=p$. Then it is
clear that $F(n;p)=\\ \mid A\mid + \mid B\mid$, where $A=\{\alpha
\in {\cal ODDP}_n: h(\alpha)=p\, \mbox{and} \, n\notin {Dom\,
\alpha}\cup {Im\, \alpha}\}$ and $B=\{\alpha \in {\cal ODDP}_n:
h(\alpha)=p\, \mbox{and} \, n\in {Dom\, \alpha}\cup {Im\,
\alpha}\}$. Define a map $\theta : \{\alpha \in {\cal ODDP}_{n-1}:
h(\alpha)=p\} \rightarrow A$ by $(\alpha)\theta=\alpha'$ where
$x\alpha'=x\alpha\,(x\in {Dom\, \alpha}$. This is clearly a
bijection since $n\notin {Dom\, \alpha}\cup {Im\, \alpha}$. Next,
recall the definitions of  $\varpi^+ (\alpha)$ and $w^+(\alpha)$
from the paragraph after Lemma \ref{lem1.1}. Now, define a map $\Phi
: \{\alpha \in {\cal ODDP}_{n-1}: h(\alpha)=p-1\} \rightarrow B$ by
$(\alpha)\Phi=\alpha'$ where

\noindent {\bf (i)} $x\alpha' =x\alpha\, (x\in {Dom\,\alpha})\,
\mbox{and}\, n\alpha' =n\,$ (if $\varpi^+ (\alpha)$ = $w^+(\alpha)$
);

\noindent {\bf (ii)} $x\alpha' =x\alpha\, (x\in {Dom\,\alpha})\,
\mbox{and}\, n\alpha' =n-\varpi^+ (\alpha)+w^+(\alpha)<n$\, (if
$\varpi^+ (\alpha) \,> w^+(\alpha)$).

\noindent In all cases $h(\alpha')=p$, and case (i) coincides with
$n\in {Dom\,\alpha'} \cap {Im\,\alpha'}$; and case (ii) coincides
with $n\in {Dom\,\alpha'} \setminus {Im\,\alpha'}$. Note that
$\varpi^+ (\alpha) \geq w^+(\alpha)$, by the order-decreasing
property. Thus $\Phi$ is onto. Moreover, it is not difficult to see
that $\Phi$ is one-to-one. Hence $\Phi$ is a bijection, as required.
This establishes the statement of the lemma. \qed

\begin{prop}\label{prop1} Let $S = {\cal ODDP}_n$ and $F(n;p)$
be as defined in (\ref{eqn1.2}) and (\ref{eqn3.1}), respectively.
Then $F(n;p)= \pmatrix{n+1\cr p+1}$, where $n\geq p \geq 1$.
\end{prop}

\pf (By Induction).

Basis Step: $F(n;1)=\pmatrix{n+1\cr 1+1}=\pmatrix{n+1\cr 2}$ and
$F(n;n)=1$ are true by Lemma \ref{lem3.1}

Inductive Step: Suppose $F(n;p)$ is true for all $n\geq p\geq 1$.

Consider $F(n+1;p)=F(n;p-1)+F(n;p)=\pmatrix{n+1\cr
p}+\pmatrix{n+1\cr p+1}\\=\pmatrix{n+2\cr p+1}=\pmatrix{(n+1)+1\cr
p+1}$, which  is the formula for $F(n+1;p)$. Hence the statement is
true for all $n\geq p\geq 1$.\qed

\begin{theorem}\label{thrm3.2} Let ${\cal ODDP}_n$ be as defined in (\ref{eqn1.2}).
Then $$\mid {\cal ODDP}_n\mid = 2^{n+1}-(n+1).$$
\end{theorem}

\pf It is enough to observe that $\mid {\cal ODDP}_n\mid=
\sum_{p=0}^{n}F(n;p)$.

\begin{lemma} \label{lem3.5} Let $S={\cal ODDP}_n$. Then
$F(n;m)={n\choose m}$, for all $n\geq m\geq 1$.
\end{lemma}

\pf It follows directly from \cite[Lemma 3.7]{Keh} and the fact that
all idempotents are necessarily order-decreasing. \qed

\begin{prop}\label{prop2} Let \,$U_n$ be a subsemigroup of ${\cal I}_n^-$ and $F(n;m)$
be as defined in (\ref{eqn3.2}). Then $F(n;0)=\mid {U_{n-1}}\mid$.
\end{prop}

\pf First, we define a map $\theta:U_{n-1} \longrightarrow
\{\alpha\in U_n: f(\alpha)=0\}$ by $\theta(\alpha)=\alpha'$ where
for all $i\,(>1)$ in ${Dom\,\alpha}$,
$$i\alpha'=(i-1)\alpha.$$
\noindent Since $n\notin {Dom\,\alpha}$ and $i\alpha'=(i-1)\alpha <
i$ for all $i>1$, it follows that $i\alpha'$ has the same degrees of
freedom as $(i-1)\alpha$, for all $i>1$. It is also clear that
$f(\alpha')=0$. Thus $\theta$ is a bijection onto $\{\alpha\in U_n:
f(\alpha)=0\}$. \qed

\begin{remark} The triangles of numbers $F(n;p)$ and $F(n;m)$,  are as
at the time of submitting this paper not in Sloane\cite{Slo}.
However,the sequence $F(n+1;m_0)$=$\mid {\cal ODDP}_n\mid$ is
\cite[A000325]{Slo}. For some computed values of $F(n;p)$ and
$F(n;m)$ in ${\cal ODDP}_n$, see Tables 3.1 and 3.2.
\end{remark}

\begin{center}

$$\begin{array}{|c|c|c|c|c|c|c|c|c|c|}
\hline
 \,\,\,\,\,n{\backslash}p&0&1&2&3&4&5&6&7&\sum F(n;p)=\mid {\cal ODDP}_n \mid
\\ \hline 0&1&&&&&&&&1
 \\ \hline 1&1&1&&&&&&&2
 \\ \hline 2&1&3&1&&&&&&5
\\
\hline 3&1&6&4&1&&&&&12
\\
\hline 4&1&10&10&5&1&&&&27
\\
\hline 5&1&15&20&15&6&1&&&58
\\
\hline 6&1&21&35&35&21&7&1&&121
\\
\hline 7&1&28&56&70&56&28&8&1&248
\\
\hline
\end{array}$$
\end{center}

\begin{center}
Table 3.1
\end{center}

\begin{center}

$$\begin{array}{|c|c|c|c|c|c|c|c|c|c|}
\hline
 \,\,\,\,\,n{\backslash}m&0&1&2&3&4&5&6&7&\sum F(n;m)=\mid {\cal ODDP}_n \mid
\\ \hline 0&1&&&&&&&&1
 \\ \hline 1&1&1&&&&&&&2
 \\ \hline 2&2&2&1&&&&&&5
\\
\hline 3&5&3&3&1&&&&&12
\\
\hline 4&12&4&6&4&1&&&&27
\\
\hline 5&27&5&10&10&5&1&&&58
\\
\hline 6&58&6&15&20&15&6&1&&121
\\
\hline 7&121&7&21&35&35&21&7&1&248
\\
\hline
\end{array}$$
\end{center}

\begin{center}
Table 3.2
\end{center}

\begin{lemma} \cite[Lemma 3.11]{Keh} Let $\alpha\in {\cal DP}_n$. Then
$\alpha$ is either order-preserving or order-reversing.
\end{lemma}

\noindent Next, we prove similar results for ${\cal DDP}_n$

\begin{lemma} \label{lem3.9} Let $\alpha\in {\cal DDP}_n$. For $1<i<n$,
if $F(\alpha)=\{i\}$ then for all $x\in {Dom\,\alpha}$ we have that
$x+x\alpha=2i$.
\end{lemma}

\pf Let $F(\alpha)=\{i\}$ and suppose $x\in {Dom\,\alpha}$.
Obviously, $i+i\alpha =i+i=2i$. If $x < i$ then $x\alpha >i$, for
otherwise we would have $i-x=\mid i\alpha-x\alpha\mid=\\ \mid
i-x\alpha\mid =i-x\alpha\Longrightarrow x=x\alpha$, which is a
contradiction. Thus, $i-x=\\ \mid i\alpha-x\alpha\mid= \mid
i-x\alpha\mid = \mid x\alpha-i\mid=x\alpha-i\Longrightarrow
x+x\alpha=2i$. The case $x>i$ is similar.\qed

\begin{lemma}\label{lem3.11} Let $S={\cal DDP}_n$. Then
$F(n;m)={n\choose m}$, for all $n\geq m\geq 2$.
\end{lemma}

\pf It follows directly from \cite[Lemma 3.18]{Keh} and the fact
that all idempotents are necessarily order-decreasing. \qed

\begin{prop}\label{prop5} Let $S={\cal DDP}_n$. Then $F(2n;m_1)=F(2n;1)=2^{n+1}-2$
and $F(2n-1;m_1)=F(2n-1;1)=3.2^{n-1}-2$, for all $n\geq 1$.
\end{prop}

\pf Let $F(\alpha)=\{i\}$. Then by Lemma \ref{lem3.9}, for any $x\in
{Dom\,\alpha}$ we have $x+x\alpha=2i$. Thus, by corollary
\ref{cor1.4}, there $2i-2$ possible elements for ${Dom\,\alpha}:
(x,x\alpha)\in \{(i,i),(i+1,i-1), (i+2,i-2), \cdots (2i-1,1)\}$.
However, (excluding $(i,i)$) we see that there are
$\sum_{j=0}{i-1\choose j}=2^{i-1}$, possible partial isometries with
$F(\alpha)=\{i\}$, where $2i-1\leq n \Longleftrightarrow i\leq
(n+1)/2$. Moreover, by symmetry we see that $F(\alpha)=\{i\}$ and
$F(\alpha)=\{n-i+1\}$ give rise to equal number of decreasing
partial isometries. Note that if $n$ is odd the equation $i=n-i+1$
has one solution. Hence, if $n=2a-1$ we have
$$2\sum_{i=1}^{a-1}2^{i-1}+2^{a-1}=2(2^{a-1}-1)+ 2^{a-1}=3.2^{n-1}-2$$
\noindent decreasing partial isometries with exactly one fixed
point; if $n=2a$ we have

$$2\sum_{i=1}^{a}2^{i-1}=2(2^a-1)=2^{a+1}-2$$
\noindent decreasing partial isometries with exactly one fixed
point.\qed

\begin{theorem}\label{thrm3.2} Let ${\cal DDP}_n$ be as defined in (\ref{eqn1.1}).
Then $$\mid {\cal DDP}_n\mid =
3a_{n-1}-2a_{n-2}-2^{\lfloor\frac{n}{2}\rfloor}+{n+1}.$$
\end{theorem}

\pf It follows from Proposition \ref{prop2}, Lemma \ref{lem3.11},
Proposition \ref{prop5} and the fact that $\mid {\cal DDP}_n\mid=
\sum_{m=0}^{n}F(n;m)$.\qed

\begin{remark} The triangles of numbers $F(n;m)$ and
the sequences $\mid {\cal DDP}_n\mid=F(n+1;m_0)$, are as at the time
of submitting this paper not in Sloane \cite{Slo}.  For some
computed values of $F(n;m)$ in ${\cal DDP}_n$, see Table 3.3.
\end{remark}

\begin{center}

$$\begin{array}{|c|c|c|c|c|c|c|c|c|c|}
\hline
 \,\,\,\,\,n{\backslash}m&0&1&2&3&4&5&6&7&\sum F(n;m)=\mid {\cal DDP}_n \mid
\\ \hline 0&1&&&&&&&&1
 \\ \hline 1&1&1&&&&&&&2
 \\ \hline 2&2&2&1&&&&&&5
\\
\hline 3&5&4&3&1&&&&&13
\\
\hline 4&13&6&6&4&1&&&&30
\\
\hline 5&30&10&10&10&5&1&&&66
\\
\hline 6&66&14&15&20&15&6&1&&137
\\
\hline 7&137&22&21&35&35&21&7&1&279
\\
\hline
\end{array}$$
\end{center}

\begin{center}
Table 3.3
\end{center}

\noindent {\bf Acknowledgements}. The first named author would like
to thank Bowen University, Iwo and Sultan Qaboos University for
their financial support and hospitality, respectively.

\vspace{1cm}

\begin{center}
R.\ Kehinde\\
Department of Mathematics and Statistics\\
Bowen University \\
P. M. B. 284,Iwo, Osun State\\
Nigeria.\\
E-mail:{\tt kennyrot2000@yahoo.com}
\end{center}

\begin{center}
S.\O.\ Makanjuola\\
Department of Mathematics\\
University of Ilorin \\
P. M. B. 1515,Ilorin, Kwara State\\
Nigeria.\\
E-mail:{\tt somakanjuola@unilorin.edu.ng}
\end{center}

\begin{center}
A.\ Umar\\
Department of Mathematics and Statistics\\
Sultan Qaboos University \\
Al-Khod, PC 123 -- OMAN\\
E-mail:{\tt aumarh@squ.edu.om}

\end{center}

\end{document}